\documentclass{amsart}
\usepackage{latexsym}
\usepackage{amsmath}
\usepackage{amssymb}
\usepackage[all]{xy}

\newtheorem{thm}{Theorem}[section]
\newtheorem{prop}[thm]{Proposition}
\newtheorem{cor}[thm]{Corollary}
\newtheorem{lem}[thm]{Lemma}
\theoremstyle{definition}
\newtheorem{defn}[thm]{Definition}

\theoremstyle{remark}
\newtheorem{rem}[thm]{Remark}

\newcommand{\K}{{\mathbb K}}
\newcommand{\F}{{\mathcal F}}
\newcommand{\e}{\varepsilon}
\newcommand{\D}{\text{D}}
\newcommand{\B}{\mathbb B}
\newcommand{\mapright}[1]{%
\smash{\mathop{%
\hbox to 1cm{\rightarrowfill}}\limits_{#1}}}
\newcommand{\maprightd}[2]{%
\smash{\mathop{%
\hbox to 1.2cm{\rightarrowfill}}\limits^{#1}\limits_{#2}}}
\newcommand{\mapleft}[1]{%
\smash{\mathop{%
\hbox to 1cm{\leftarrowfill}}\limits_{#1}}}
\newcommand{\mapleftu}[1]{%
\smash{\mathop{%
\hbox to 0.8cm{\leftarrowfill}}\limits^{#1}}}
\newcommand{\maprightu}[1]{%
\smash{\mathop{%
\hbox to 1cm{\rightarrowfill}}\limits^{#1}}}
\newcommand{\maprightud}[2]{%
\smash{\mathop{%
\hbox to 1cm{\rightarrowfill}}\limits^{#1}_{#2}}}
\newcommand{\mapleftud}[2]{%
\smash{\mathop{%
\hbox to 1cm{\leftarrowfill}}\limits^{#1}_{#2}}}

\newcounter{eqn}[section]

\def\theeqn{\textnormal{(\thesection.\arabic{eqn})}}

\def\eqnlabel#1{%
\refstepcounter{eqn}%
\label{#1}%
\leqno{\theeqn}}

\begin{document}

\title[The Hochschild cohomology ring of a singular cochain algebra]
{The Hochschild cohomology ring of the singular cochain 
algebra of a space}

\footnote[0]{ 


Department of Mathematical Sciences,
Faculty of Science,
Shinshu University,
Matsumoto, Nagano 390-8621, Japan
e-mail:{\tt kuri@math.shinshu-u.ac.jp}
}

\author{Katsuhiko KURIBAYASHI}
\date{}

\begin{abstract} 
We determine the algebra structure of the Hochschild cohomology 
of the singular cochain algebra with coefficients in a field  
on a space whose cohomology is a polynomial algebra.
A spectral sequence calculation of the Hochschild cohomology
is also described. In particular, when the underlying field is 
of characteristic two,
we determine the associated bigraded Batalin-Vilkovisky algebra structure on 
the Hochschild cohomology of the singular cochain on a space 
whose cohomology is an exterior algebra. 
\end{abstract}

\maketitle

\section{Introduction}

The objective of this paper is to compute the Hochschild cohomology ring
of the singular cochain algebra on a simply-connected space,
whose cohomology is either a polynomial algebra or an exterior algebra.
The motivational topics are in string topology as well as 
in the classification problem 
of thick subcategories of the triangulated category
associated with a singular cochain algebra on a space.

Let $M$ be a compact, oriented $d$-dimensional smooth manifold and 
$LM=map(S^1, M)$ the space of free loops on $M$.
In \cite{C-S}, Chas and Sullivan have defined
a product on the shifted homology
${\mathbb H}_*(LM):=H_{*+d}(LM)$, which is called the loop homology of $M$,
and have shown that the homology, moreover, has the structure
of a Batalin-Vilkovisky algebra.
A result due to Cohen and Jones in \cite{C-J} asserts that there exists 
an isomorphism of algebras between the loop homology of $M$ and
the Hochschild cohomology ring of the singular cochain algebra on $M$;
see also \cite{C}. 
This allows one to describe the loop homology 
in terms of differential homological algebra. For various developments along
this line, we refer the reader to \cite{F-T1, F-T2, F-M-T, F-T-V,
F-T-V2, Kauf, Kauf_2, Luc1, Luc2, Me, T, T-Z}. 
Especially, 
Menichi \cite{Luc2} has shown that 
the Hochschild cohomology of the singular cochain on a Poincar\'e
duality space admits the structure 
of a Batalin-Vilkovisky algebra with the B-V operator given rise to by
the Connes coboundary map.  
Thus we are also led to the study of the extended structure 
of the Hochschild cohomology. 

As for global nature of singular cochains, 
J{\o}rgensen has investigated 
the derived category of
the singular cochain on a Poincar\'e duality space
by applying Auslander-Reiten theory. In particular,
the Auslander-Reiten quiver of the full subcategory consisting of
compact objects is determined 
in \cite{Jo} and \cite{Jo2}, see also \cite{Schmidt}.
Such the result brings us to the study of topological spaces
with categorical representation theory. 
Indeed, it is expected that numerical invariants, 
which appear and grow uniquely in the representation theory, capture
topological properties of spaces via functors from the category of spaces
to an algebraic one. 
The author has introduced 
in \cite{K2} a new topological invariant, 
which is called {\it the level of a space}. For a space $M$ 
over a given space, the level of $M$ measures the number of steps
to build the singular cochain on the space $M$ from that 
on a more fundamental space under an appropriate rule. 
We refer the reader to \cite{K3} 
for a linkage  between the level and
the Lusternik-Schnirelmann category of a space. 

The origin of the numerical invariant level is in 
the study of the dimensions of triangulated categories; 
see for example \cite{R}. Then {\it the level of an object}
in a triangulated category was first introduced
by Avramov, Buchweitz, Iyengar and Miller in \cite{ABIM}. 
It is also mentioned that, relying on knowledge of the levels of
vertices on the Auslander-Reiten quiver alluded to above, the explicit
calculation of the level for the total space of a bundle over the sphere
is performed in \cite{K2}. 

Recently, Benson, Iyengar and Krause \cite{B-I-K} have proved
a classification theorem of thick subcategories of
a triangulated category endowed with a ring homomorphism
from a (graded) commutative ring
to the graded center of the category; see Appendix for terminology. 
Let $C^*(X;\K)$ be the singular cochain algebra on a space $X$ 
with coefficients in a field $\K$ and $\D(C^*(X;\K))$ the derived
category of differential graded modules over 
the differential graded algebra $C^*(X;\K)$.
Then the level of a space over $X$
is defined in the category $\D(C^*(X;\K))$, 
which has the structure of a triangulated category. 
Therefore, we also expect that an explicit classification of thick
subcategories of $\D(C^*(X;\K))$ plays an important role in the study of
the levels of spaces because the invariants   
are defined by filtering a suitable thick subcategory of $\D(C^*(X;\K))$. 
In particular, one can take the Hochschild cohomology ring of $C^*(X; \K)$
as the graded commutative ring that may be the key to the
classification of thick subcategories of 
the triangulated category $\D(C^*(X;\K))$.  

Thus the Hochschild cohomology of the singular
cochain on a space becomes one of our great interests.
However, there are few results on explicit calculations
of the Hochschild cohomology rings
of singular cochain algebras $C^*(X;\K)$ except for the case that
the cohomology algebra $H^*(X; \K)$ is generated by a single element;
see \cite{Luc1, Yang}.

As mentioned above, in this paper, we confine our attention mainly to the
calculation of the Hochschild cohomology $HH^*(C^*(X; \K);C^*(X;\K))$ 
of the singular cochain on a
simply-connected space $X$ whose cohomology is either a polynomial algebra
or an exterior algebra. Unless otherwise explicitly stated, it is
assumed that a space has the homotopy type of a CW-complex. 
%

Our results in this paper are now described. 

\begin{thm}
\label{thm:1}
Let $X$ be a simply-connected space whose mod $p$ cohomology is
a polynomial algebra, say
$
H^*(X; {\mathbb Z}/p) \cong {\mathbb Z}/p[x_1, x_2, ..., x_n].
$
Then
$$
HH^*(C^*(X;{\mathbb Z}/p); C^*(X;{\mathbb Z}/p))
\cong {\mathbb Z}/p[x_1, x_2, ..., x_n]\otimes
\wedge (u_1^*, u_2^*, ..., u_n^*)
$$
as an algebra, where $\deg u_i ^*= - \deg x_i +1$.
\end{thm}

Before describing another result on the computation 
of the Hochschild cohomology, we here recall the definition of
the Batalin-Vilkovisky algebra.

\begin{defn}
A commutative graded algebra $A^*$ is a Batalin-Vilkovisky algebra 
if $A^*$ is equipped with an operation 
$\Delta : A^* \to A^{*-1}$ such that $\Delta^2=0$ and 
\begin{eqnarray*}
\Delta(abc) &=& \Delta(ab)c + (-1)^{|a|}a\Delta(bc) +  
(-1)^{(|a|-1)|b|}b\Delta(ac) \\
                & & \ \ \  -(\Delta a)bc - (-1)^{|a|}a(\Delta b) c 
- (-1)^{|a|+|b|}ab(\Delta c). 
\end{eqnarray*}
The map $\Delta$ is called the B-V operator. 
\end{defn}

Observe that the action of the B-V operator on the product of three elements is
determined exactly by knowledge of that on the product of two elements.    

Let $\K$ be a fixed field. A path-connected space $M$ is called 
{\it a Poincar\'e duality space} of formal dimension $m$ if
the space $M$ is equipped with an {\it orientation class}
$[M] \in H_m(M; \K)$ such that the cap product
$$
- \cap [M] : H^*(M; \K) \to H_{m-*}(M; \K)
$$
is an isomorphism. {\it The fundamental class} of $M$ is the element 
$\omega_M$ such that $\langle \omega_M, [M] \rangle= 1$, where 
$\langle \ , \ \rangle$ denotes the Kronecker product. 

In what follows, for a space $X$, we may write $C^*(X)$ for $C^*(X; \K)$.
Let $M$ be a simply-connected Poincar\'e duality space. 
Theorem \ref{thm:MSS} below states that 
the Moore spectral sequence (MSS for short) converges to   
$HH^*(C^*(M); C^*(M))$ as a Batalin-Vilkovisky algebra.    
More precisely, each term of the MSS admits a
differential Batalin-Vilkovisky algebra structure and 
the induced filtration on $HH^*(C^*(M); C^*(M))$ respects the B-V operator 
defined by Menichi \cite{Luc2}. Moreover, the $E_\infty$-term is
isomorphic to the bigraded algebra 
$\text{Gr}HH^*(C^*(M); C^*(M))$ associated with the filtration 
as a bigraded Batalin-Vilkovisky algebra equipped with   
the B-V operator of bidegree $(-1, 0)$. 

The MSS with the differential Batalin-Vilkovisky algebra structure
provides a new method for computing the Hochschild cohomology 
of the singular cochain on a space.  
In fact by applying the spectral sequence, we establish  

\begin{thm}
\label{thm:2}Let $M$ be a simply-connected space whose 
mod $2$ cohomology is an exterior algebra, say
$
H^*(M; {\mathbb Z}/2) \cong \wedge(y_1, y_2, ..., y_l).
$
Suppose further that the operation $Sq^1$ vanishes on the cohomology. 
Then as a bigraded Batalin-Vilkovisky algebra, 
$$
\text{\em Gr}HH^*(C^*(M;{\mathbb Z}/2); C^*(M;{\mathbb Z}/2)) \cong
\wedge(y_1, y_2, ..., y_l)\otimes {\mathbb Z}/2[\nu_1^*, \nu_2^*, ...,
\nu_l^*]
$$
in which $\Delta(y_j)=0$, $\Delta(\nu_i^*)=0$, $\Delta(y_iy_j)=0$, 
$\Delta(\nu_i^*\nu_j^*)=0$ 
for $1\leq i, j \leq l$ and 
$\Delta(y_i\nu_j^*)=\delta_{ij}\cdot 1$, where 
$\text{\em bideg} \  y_j =(0, \deg y_j)$ and $\text{\em bideg} \
\nu_j^* 
=(1, -\deg y_j)$ for $1\leq j  \leq l$.
\end{thm}

Observe that the space $M$ as in Theorem \ref{thm:2} is 
a Poincar\'e duality space with a orientation class $[M]$ which
is the dual  to a top non-zero element of the cohomology. 

For a very special case, we can solve the extension problems 
on the product and on the B-V operator, which appear in the bigraded 
Batalin-Vilkovisky algebra 
$\text{Gr}HH^*(C^*(M;{\mathbb Z}/2); C^*(M;{\mathbb Z}/2))$; 
see Corollary \ref{cor:S-times-S}. 
It seems that the result is the first computational example obtained by
means of the MSS. 

\begin{rem}In general, the squaring operation $Sq^1$ acts non-trivially
on the mod $2$ cohomology algebra of the Stiefel manifold of 
the form $M=SO(m+n)/SO(n)$ even if the algebra
is an exterior algebra, see \cite{M-T}.
However, thanks to the results \cite[Corollary 5]{K1} and 
Proposition \ref{prop:key} below,
we see that the conclusion of Theorem \ref{thm:2} remains valid for 
the Stiefel manifold $SO(m+n)/SO(n)$ provided $m\leq \min \{4, n\}$.
\end{rem}

The rest of this paper is organized as follows.
In Section 2, we recall the cup product of the Hochschild cohomology of
a differential graded algebra and prove Theorem \ref{thm:1}.
Section 3 is devoted to proving the assertion in Theorem \ref{thm:2}
concerning the bigraded algebra structure. To this end, we review the
Moore spectral sequence converging to the Hochschild cohomology of the
singular cochain algebra on a space. 
In Section 4, we discuss a Batalin-Vilkovisky algebra structure on
the spectral sequence. Moreover fundamental properties  
of the Moore spectral sequence are considered here. 
In consequence, Theorem \ref{thm:2} is proved completely.
Though Appendix has no result,  
we describe how one can take the loop homology into the categorical
representation theory via the Hochschild cohomology.

\section{The cup product in the Hochschild cohomology}

We begin with the definition of the Hochschild cochain complex. 
Let $(A, d)$ be an augmented differential graded algebra over a field $\K$ and
$s\bar{A}$ denote the suspension of the augmentation ideal $\bar{A}$; that is, 
$(s\bar{A})^n =\bar{A}^{n+1}$.
Let $T(s\bar{A})$ stand for the tensor algebra on $s\bar{A}$.
The two-sided normalized bar complex ${\mathbb B}(A;A;A)$ is the tensor product
$A\otimes T(s\bar{A}) \otimes A$
with the differential $d_{{\mathbb B}}= d_1+d_2$ defined by
\begin{eqnarray*}
d_1(a[a_1|a_2|...|a_k]b) &=& d(a)[a_1|a_2|...|a_k]b
- \sum_{i=1}^{k}(-1)^{\e_i}
a[a_1|a_2|...|d(a_i)|...|a_k]b \\
& & +(-1)^{\e_{k+1}}a[a_1|a_2|...|a_k]d(b), \\
d_2(a[a_1|a_2|...|a_k]b) &=& (-1)^{|a|}aa_1[a_2|...|a_k]b
+ \sum_{i=2}^{k}(-1)^{\e_i}a[a_1|a_2|...|a_{i-1}a_i|...|a_k]b \\
& & -(-1)^{\e_{k}}a[a_1|a_2|...|a_{k-1}]a_kb,
\end{eqnarray*}
where $\e_i= |a| + \sum_{j<i}(|sa_j|)$.

Let $(N, d_N)$ be a differential graded $A$-module.
Then by definition the Hochschild cochain complex is a complex
${\bf C}(A, N)=\{{\bf C}^n(A, N), \partial\}$
of the form
$$
{\bf C}^n(A, N)=\text{Hom}_{A\otimes A^{op}}^n({\mathbb B}(A;A;A), N)
$$
with the differential $\partial$  defined by
$\partial(f)=d_N f - (-1)^{|f|}f d_{{\mathbb B}}$. The Hochschild cohomology
$HH^*(A; N)$ is defined
to be the homology of the complex ${\bf C}(A,N)$.
It follows from \cite[Proposition 19.2]{F-H-T} that the multiplication
on $A$
induces a semi-free resolution
${\mathbb B}(A;A;A) \stackrel{\simeq}{\to} A$ of $A\otimes
A^{op}$-modules. 
This enables us to regard the Hochschild cohomology as the Ext-group
$\text{Ext}_{A\otimes A^{op}}(A, N)$ in the sense of Moore; see
\cite[Appendix]{F-H-T_G}.

For a vector space $V$, 
we denote by $V^\vee$ the dual vector space to $V$ 
unless otherwise noted.   
For a basis $\{v_i\}_{i\in I}$ for $V$, 
the dual basis is denoted by  $\{v_i^*\}_{i\in I}$. 

We here recall a Koszul-Tate resolution of a graded commutative algebra.
Let
$\Lambda$ a graded algebra over a field $\K$ of the form
$$
\Lambda = \wedge (y_1, ..., y_l)\otimes
\K[x_1, ..., x_n]/(\rho_1, ..., \rho_m),
$$
where $\rho_1, ..., \rho_m$ is a regular sequence in
the polynomial algebra $\K[x_1, ..., x_n]$.
We assume that each $\rho_i$ is decomposable.
Observe that  $\deg x_i$
is even and $\deg y_i$ is odd 
if the characteristic of $\K$ is greater than two. 
The algebra $\Lambda$ is called
{\it a graded complete intersection algebra}.

\begin{prop}{\em (}\cite[Proposition 3.5]{Smith1}
\cite[Proposition 1.1]{K1} {\em )}
\label{prop:resolution}
Under the above hypothesis, there exists a projective
resolution ${\mathcal F} \stackrel{\varphi}{\to}
\Lambda \to 0$ of $\Lambda$ as a left
$\Lambda\otimes \Lambda^{op}$-module such that
$$
{\mathcal F} = \Lambda\otimes \Lambda \otimes
\Gamma[\nu_1, ..., \nu_l]\otimes \wedge (u_1, ..., u_n) \otimes
\Gamma[w_1, ..., w_m],
$$
$d(\Lambda\otimes \Lambda)=0$, $d(\nu_i)=y_i\otimes 1 - 1\otimes y_i$,
$d(u_j)=x_j\otimes 1 -1 \otimes x_j$,
$d(\gamma_r(w_i))=(\sum_{j=1}^n\zeta_{ij}u_j)\otimes \gamma_{r-1}(w_i)$
and $\varphi$ is the multiplication of $\Lambda$, where 
$\text{\em bideg} \ \lambda = (0, \deg \lambda)$ for
$\lambda \in \Lambda\otimes \Lambda$,
$\text{\em bideg} \ \nu_i = (-1, \deg y_i)$,
$\text{\em bideg} \ u_j = (-1, \deg x_i)$ and
$\text{\em bideg} \ w_i = (-2, \deg \rho_i)$. Here
$\zeta_{ij}$ is an appropriate element of
$\K[x_1, ..., x_n]\otimes \K[x_1, ..., x_n]$
which satisfies the condition that
$$
\rho_i\otimes 1 - 1\otimes \rho_i =
\sum_{j=1}^n\zeta_{ij}(x_j\otimes 1- 1\otimes x_j) \ \
\text{and} \ \ \varphi(\zeta_{ij})= \frac{\partial \rho_i}{\partial
x_j}.
$$
\end{prop}

It is readily seen that ${\mathcal F}$ is semi-free and hence 
$\varphi : {\mathcal F} \to \Lambda$ in
Proposition \ref{prop:resolution} is a semi-free resolution of
$\Lambda$.  
In what follows, we shall call the resolution 
the Koszul-Tate resolution of $\Lambda$.

\begin{rem}
In the case where the algebra $\Lambda$ is
a truncated polynomial algebra generated by a single element, the
Koszul-Tate resolution is nothing but the periodic resolution, which is
used in \cite{Yang, Holm, Ci-So} to compute the Hochschild
(co)homology of $\Lambda$. See also \cite{H-X} for the Hochschild
cohomology ring of an exterior algebra.
\end{rem}

Let $A$ be an augmented differential graded algebra (DG algebra), 
$M$ and $N$ differential graded bimodules over $A$. 
Let $P \stackrel{\simeq}{\to} A$ be a
semi-free resolution of $A$ as a $A\otimes A^{op}$-module. 
Then the cup product
$$
\smile \ : HH^*(A; M)\otimes HH^*(A; N) \to HH^*(A, M\otimes_AN)
\eqnlabel{add-0}
$$
is defined with an $A\otimes A^{op}$-chain map
$
D  : P \to P\otimes_AP,
$
which is a lift of the identity map on $A$, by the composite
$$
f\smile g : P \stackrel{D }{\to} P\otimes_A P
\stackrel{f\otimes g}{\to} N\otimes_A N
$$
for $f, g \in \text{Hom}_{A\otimes A^{op}}(A, N)$.
We call the map $D  : P \to P\otimes_AP$ {\it a diagonal map}.

Suppose that there exists an $A\otimes A^{op}$-module map
$N\otimes _A N \to N$. Then the cup product on
the Hochschild cochain complex ${\bf C}(A; N)$ makes
the Hochschild cohomology $HH^*(A; N)$ into a graded algebra.
Observe that, in the case $M=N=A$, the algebra 
$HH^*(A; A)$ is graded commutative; see for example 
\cite[Proposition 1.2]{Sanada}. Moreover, using the $A\otimes A^{op}$-module isomorphism  
$A\otimes_AN \stackrel{\cong}{\to} N$, which is induced by the right $A$-module structure on $N$, we give 
$HH^*(A; N)$ an $HH^*(A; N)$-module structure 
$$
\smile \ : HH^*(A; A)\otimes HH^*(A; N) \to HH^*(A, N). 
\eqnlabel{add-1}
$$

We here look at a bigraded algebra structure on the Hochschild
cohomology
$HH^*(A; N)$ provided the differentials on $A$ and $N$ are trivial.
Let ${\bf C}^{p,q}$ be the subspace
$\text{Hom}_{A\otimes A^{op}}^{p+q}(P_{-p}, N)$ of ${\bf C}(A; N)$.
Then it follows that the differential $d$ of ${\bf C}(A; N)$
maps ${\bf C}^{p,q}$ into ${\bf C}^{p+1,q}$.
This implies that $\{ {\bf C}^{*,q}, d\}_{p\geq 0}$
is a subcomplex of the Hochschild cochain complex
${\bf C}(A; N)$ for any $q$ and hence we have
$$
HH^n(A;N) =\bigoplus_{p+q=n}HH^{p,q}(A, N),
$$
where $HH^{*,q}(A, N)$ denotes the homology of the complex $\{ {\bf
C}^{*,q},
d\}_{p\geq 0}$.
Let $D  : P \to P\otimes_A P$ be a diagonal map. 
Then by definition one sees that the image  
$D(P _{-p})$ is included in $\oplus_{i+j=-p} P_i\otimes_A P_j$.
This yields that the Hochschild cohomology $HH^n(A;N)$ admits
a bigraded algebra structure; that is,
$$
\smile \ : HH^{p,q}(A, N)\otimes HH^{p',q'}(A, N) \to HH^{p+p',q+q'}(A, N).
$$

The uniqueness of the cup product on the Hochschild homology; see
\cite[\S1]{Sanada}\cite[\S2]{S-W},
allows us to define the product with a tractable diagonal map.
We construct an explicit diagonal map for
the Koszul-Tate resolution of a DG algebra, 
which is the tensor product of a polynomial
algebra and an exterior algebra with the trivial differential.

Suppose that $\Lambda$ is an algebra of the form
$
\wedge (y_1, ..., y_l)\otimes \K[x_1, ..., x_n].
$
We write $\Lambda\otimes \Lambda \otimes {\mathcal E}$
for the graded algebra ${\mathcal F}$ 
mentioned in Proposition \ref{prop:resolution} with
${\mathcal E}=\Gamma[\nu_1, ..., \nu_l]\otimes \Lambda (u, ..., u_n)$.
Define a $\Lambda\otimes \Lambda^{op}$-homomorphism 
$$D  : \Lambda\otimes \Lambda \otimes {\mathcal E}
\to \Lambda\otimes \Lambda \otimes {\mathcal E}\otimes_\Lambda
\Lambda\otimes \Lambda \otimes {\mathcal E}, 
$$ 
on the generators of the algebra 
${\mathcal F}=\Lambda\otimes \Lambda \otimes {\mathcal E}$ 
by 
$$
D (u_j)= 1\otimes 1 \otimes u_j\otimes_\Lambda 1\otimes 1\otimes 1 +
1\otimes 1\otimes 1\otimes_\Lambda 1\otimes 1\otimes u_j,
$$
$$
D (\gamma_k(\nu_i)) =\sum_{s+t=k}
1\otimes 1\otimes \gamma_s(\nu_i)\otimes_\Lambda
1\otimes 1\otimes \gamma_t(\nu_i)
$$ 
and extend them to the whole module with 
$$
D (\alpha \gamma_{i_1}(\nu_1)\cdots \gamma_{i_l}(\nu_l)u_1^{\e_1}\cdots u_n^{\e_n})=
\alpha D(\gamma_{i_1}(\nu_1))\cdots D(\gamma_{i_l}(\nu_l))
D(u_1)^{\e_1}\cdots D(u_n)^{\e_n}, 
$$
where $\e_i= 0$ or $1$ and $\alpha \in \Lambda\otimes \Lambda^{op}$. 
Observe that 
$$
\lambda_1\otimes \lambda_2 (\lambda \otimes \mu \otimes a \otimes_\Lambda
\lambda' \otimes \mu'\otimes b)
=(-1)^{|\lambda_2|(|\lambda|+|\mu|+|a|+|\lambda'|)}
\lambda_1\lambda \otimes \mu \otimes a \otimes_\Lambda  
\lambda' \otimes \lambda_2 \mu' \otimes b, 
$$
for $\lambda_1\otimes \lambda_2 \in \Lambda\otimes \Lambda^{op}$ and 
$\lambda \otimes \mu \otimes a \otimes_\Lambda  
\lambda' \otimes \mu'\otimes b \in {\mathcal F}\otimes_\Lambda {\mathcal F}$.

\begin{lem}
The $\Lambda\otimes \Lambda^{op}$-homomorphism $D $ is a diagonal map.
\end{lem}

\begin{proof}
The differential $d$ is a derivation on the algebra ${\mathcal F}$ and
hence so is $\partial:=d\otimes 1 + 1\otimes d$ on
${\mathcal F}\otimes_\Lambda {\mathcal F}$.
Thus in order to prove the lemma, it suffices to show
that 
$$\partial D(v) = D d(v)
\eqnlabel{add-2}
$$ for any
$v \in \{\gamma_{i_1}(\nu_1), ..., \gamma_{i_l}(\nu_l)\ | i_t \geq 1\} 
\cup \{u_1, ..., u_n, x_i \}$. In fact, we choose a base $x$ for ${\mathcal F}$ of the form  
$\alpha v_i\cdots v_l\cdot v_{l+1}\cdots v_{l+k}$, where 
$\alpha \in   \Lambda\otimes \Lambda^{op}$, 
$v_j=\gamma_{i_j}(\nu_j)$, $v_{l+i} \in  \{u_1, ...., u_n\}$ and $v_{l+i}\neq v_{l+j}$ if $i\neq j$. Then we see that 
\begin{eqnarray*}
\partial D (x) &=& \partial (\alpha D(v_1)\cdots D(v_{l+k})) \\
&=& \alpha \sum_{1\leq i \leq l+k} \pm D(v_1) \cdots \partial D(v_i)\cdots D(v_{l+k})\\
&=& \alpha \sum_{1\leq i \leq l+k} \pm D(v_1) \cdots D(d v_i)\cdots D(v_{l+k})\\
&=& D(\alpha \sum_{1\leq i \leq l+k} \pm v_1 \cdots d v_i \cdots  v_{l+k}) \ = \ Dd(x).
\end{eqnarray*} 
Here $\pm$ denotes the Koszul sign. 
The differential $d$ is closed under the subalgebras 
$\Lambda\otimes \Lambda\otimes \Gamma[\nu_i]$ and $\Lambda\otimes \Lambda\otimes \wedge (u_j)$ 
for $i = 1, .., l$ and $j=1, ..., n$ so that 
the forth equality follows from the definition of the diagonal map $D$. We now verify the equality (2.3). 
It follows that
\begin{eqnarray*}
& & \partial D (\gamma_k(\nu_i)) = \partial
(\sum_{s+t=k}1\otimes 1 \otimes \gamma_s(\nu_i)\otimes _\Lambda 1\otimes
1\otimes \gamma_t(\nu_i)) \\
&=&
\sum_{s+t=k}\big(y_i\otimes 1 \otimes \gamma_{s-1}(\nu_i)\otimes _\Lambda
1\otimes 1\otimes \gamma_t(\nu_i)
- 1\otimes y_i \otimes \gamma_{s-1}(\nu_i)\otimes _\Lambda \otimes 1
\otimes \gamma_t(\nu_i) \\
& & + 1\otimes 1 \otimes \gamma_s(\nu_i)\otimes _\Lambda y_i\otimes 1
\otimes \gamma_{t-1}(\nu_i)
- 1\otimes 1 \otimes \gamma_s(\nu_i)\otimes _\Lambda 1\otimes y_i
\otimes \gamma_{t-1}(\nu_i) \big)\\
&=& \sum_{s+t=k-1}\big(y_i\otimes 1 \otimes \gamma_{s}(\nu_i)\otimes
_\Lambda 1\otimes 1\otimes \gamma_t(\nu_i)
- 1\otimes 1 \otimes \gamma_s(\nu_i)\otimes _\Lambda 1\otimes y_i
\otimes \gamma_{t}(\nu_i) \big) \\
&=& (y_i\otimes 1 - 1\otimes y_i) \cdot \big(
\sum_{s+t=k-1}1\otimes 1 \otimes \gamma_{s}(\nu_i)\otimes _\Lambda
1\otimes 1\otimes \gamma_t(\nu_i)
\big) \\
&=& D  ((y_i\otimes 1 - 1\otimes y_i)\gamma_{k-1}(\nu_i)) = D 
d(\gamma_k(\nu_i))
\end{eqnarray*}
The same calculation as above enables us to conclude  that
$\partial D (u_j)= D  d (u_j)$.
We have the result.
\end{proof}

Let $A$ be a $\Lambda$-bimodule equipped with 
a  $\Lambda \otimes \Lambda^{op}$-module map $A\otimes_\Lambda A \to A$.  
In particular, if $A$ is a commutative graded algebra over $\Lambda$,
then $A$ is viewed as a $\Lambda$-bimodule with 
$(\lambda_1 a)\lambda_2 =\lambda_1 (a \lambda_2) 
:= (-1)^{|a||\lambda_2|}\lambda_1(\lambda_2 a)$  
for $\lambda_1, \lambda_2 \in \Lambda$ and $a\in A$. 
Moreover, a $\Lambda \otimes \Lambda^{op}$-module map 
$A\otimes_\Lambda A \to A$ is naturally induced by the product on $A$. 

Assume that
$\Lambda$ is isomorphic to either a polynomial algebra 
$\K[x_1, ..., x_n]$ or an exterior algebra
$\wedge(y_1, ..., y_l)$ and that $A$ is a commutative graded algebra over $\Lambda$ such that 
$\dim A < \infty$ if $\Lambda$ is an exterior algebra.
We consider the DG algebra structure of the complex $\text{Hom}_{\Lambda
\otimes \Lambda}({\mathcal F}, A)$,
which computes the $E_2$-term of the spectral sequence 
introduced in the next section.

By assumption, one of algebras $A$ and ${\mathcal E}$ is of finite dimension. 
This allows us to obtain an isomorphism 
$$
\theta : A\otimes {\mathcal E}^\vee = A\otimes \text{Hom}_\K({\mathcal
E}, \K) \stackrel{\cong}{\to}
\text{Hom}_{\Lambda \otimes \Lambda}(\Lambda \otimes \Lambda\otimes
{\mathcal E}, A)
$$
defined by $\theta(a\otimes f)(\lambda \otimes \alpha)
=(-1)^{|\lambda|(|f| + |a|)}\lambda \cdot a
f(\alpha)$, where $a \in A$,
$\lambda \in \Lambda \otimes \Lambda$,
$f \in {\mathcal E}^\vee$, $\alpha \in {\mathcal E}$ and
$\cdot$ stands for
the $\Lambda \otimes \Lambda$-module structure on $A$.
Thus it follows that the vector space $A\otimes {\mathcal E}^\vee$ admits a
differential algebra structure via
the isomorphism $\theta$.
As for the algebra structure, we see that for the dual base
$\gamma_k(\nu_i)^*$ to $\gamma_k(\nu_i)$ and
the dual base $u_j^*$ to $u_j$,
$$
a\otimes \gamma_k(\nu_i)^* \cdot b \otimes \gamma_l(\nu_i)^* = ab \otimes
\gamma_{k+l}(\nu_i)^*
\ \ \text{and} \ \
a\otimes u_j^* \cdot b\otimes u_j^* =0.
$$
Observe that the total degree of the element $\gamma_k(\nu_i)^*$ is even
if $p$ is odd.  
For any element $x$ in $\Lambda\otimes \Lambda \otimes {\mathcal E}$, each term of $dx$ 
has an element of the form $ \lambda \otimes 1 - 1\otimes \lambda$, 
where $\lambda \in \Lambda$; see  Proposition \ref{prop:resolution}. 
Moreover, since 
$(\lambda \otimes 1 - 1\otimes \lambda)\cdot a = \lambda (1a) -
1(\lambda a) = 0$ 
for $\lambda \in \Lambda$ and $a \in A$, it follows that for any 
$\Lambda\otimes \Lambda$-module map 
$\varphi : \Lambda\otimes \Lambda \otimes {\mathcal E} \to A$, 
$(d\varphi)(x)=(-1)^{|\varphi |}\varphi(dx)=0$ 
and hence the differential on $A\otimes
{\mathcal E}^\vee$ is trivial.
Thus we have 

\begin{prop}
\label{prop:HH}
As a bigraded algebra,
$$
HH^*(\Lambda; A) \cong
\left\{
\begin{array}{ll}
A\otimes \K[\nu_1^*, ...., \nu_l^*] & \ \text{if} \ \ \Lambda =
\wedge(y_1, .., y_l), \\
A \otimes \wedge (u_1^*, ..., u_n^*) & \ \text{if} \ \ \Lambda = \K[x_1,
.., x_n], 
\end{array}
\right. 
$$
where $\text{\em bideg} \ a = (0, \deg a)$, $\text{\em bideg} \ \nu_i^* = (1,
-\deg y_i)$ and
$\text{\em bideg} \ u_j^* = (1, -\deg x_j)$.
\end{prop}

\noindent
{\it Proof of Theorem \ref{thm:1}.}
The proof of \cite[7.1 Theorem]{Mun}
implies that $C^*(X)$ is $\K$-formal; that is, there exists a sequence
of quasi-isomorphisms which connects $C^*(X)$ with $H^*(X;\K)$. It follows
from \cite[3.4 Proposition]{F-M-T} that the Hochschild cohomology ring
$HH(C^*(X); C^*(X))$ is isomorphic to
$HH(H^*(X); H^*(X))$ as an algebra.
Proposition \ref{prop:HH} yields the result. \hfill \qed

\begin{rem}
Let $\Lambda$ be a graded complete intersection algebra.
The algebra structure of  the Hochschild cohomology $HH^*(\Lambda, \Lambda)$ may be
described in terms of  cycles on the Koszul-Tate resolution in
Proposition \ref{prop:resolution} if one has an explicit form of a
diagonal map.
\end{rem}

\section{A spectral sequence converging to
the Hochschild cohomology ring of a DG algebra}

For a space $X$, we assume that the cohomology $H^*(X; \K)$ is of finite
type; that is, $\dim H^i(X; \K) < \infty$ for any $i$.    
Let $M$ and $N$ be connected spaces 
and $f : N \to M$  a map. The singular cochain algebra
$C^*(N)$ is regarded as a $C^*(M)$-bimodule via the map
$f^* : C^*(M) \to C^*(N)$ induced by $f$.
Then it follows that the cup product gives rise to a 
$C^*(M)\otimes C^*(M)^{op}$-module map 
$C^*(N)\otimes_{C^*(M)}C^*(N) \to C^*(N)$.  

\begin{thm} \text{\em (}cf. \cite[1 Proposition]{F-T-V} \text{\em )}
\label{thm:Mss} 
Under the above hypothesis, we assume further that
$H^*(N)$ is of finite dimension.
Then there exists a right-half plane cohomological
spectral sequence $\{E_r^{*,*}, d_r\}$ converging to the Hochschild
cohomology $HH^*(C^*(M); C^*(N))$ as an algebra such that
$$
E_2^{p,q} \cong HH^{p, q}(H^*(M); H^*(N))
$$
as a bigraded algebra.
\end{thm}

\begin{proof} 
Let $S$ be a complement of the vector subspace generated by cycles of 
$C^d(N)$, where 
$d = \sup \{n \ | \ H^*(N)\neq 0 \}$. 
We define $I$ to be the two-sided ideal generated by $C^{>d}(N)\oplus S$. 
Then the projection $C^*(N) \to C^*(N)/I$ is 
a quasi-isomorphism of $A$-bimodules. 

Let $A:=TV \stackrel{\simeq}{\to} C^*(M)$ be a TV model for the space $M$ 
in the sense of Halperin and Lemaire \cite{H-L}. Let $\B_*(A; A; A)$ be 
the normalized bar complex mentioned in the previous section.  
We then define
a decreasing filtration $\{F^p{\bf C}^*\}_{p \geq 0}$
of the Hochschild cochain complex
$
{\bf C}^*=\{\text{Hom}_{A\otimes A^{op}}
(\B_*(A;A;A), C^*(N)/I) \}_{n \in {\mathbb Z}}
$ by
$$
F^p{\bf C}^n =
\prod_{s\geq p} \text{Hom}_{A\otimes A^{op}}^n(\B_s(A;A;A), C^*(N)/I), 
$$
where $\B_s(A;A;A)=A\otimes s\overline{A}^{\otimes s}\otimes A$.
Since $s\overline{A}$ has no element of degree zero, it follows that the number 
$\inf \{ m \ | \ (\B_s(A;A;A))^m \neq 0\}$ increases strictly if so does $s$. It is immediate that 
$(C^*(N)/I)^{< 0} = 0$ and $(C^*(N)/I)^{> d}=0$. 
These facts imply that 
the filtration $\{F^p{\bf C}^*\}_{p \geq 0}$ is bounded; that is, for any $n$, there exists $p(n)$ such that 
$F^{p}{\bf C}^n = 0$ for $p > p(n)$.

Observe that $F^0{\bf C}^* = {\bf C}^*$ and that
the cup product on the Hochschild cochain complex
respects the filtration; that is,
$F^s{\bf C}^n \smile F^t{\bf C}^m \subset F^{s+t}{\bf C}^{n+m}$.
Therefore as usual 
we can construct a spectral sequence $\{E_r, d_r\}$, whose each term
admits the structure of a differential graded algebra, 
by using the filtration.  Since the filtration is bounded, it follows that the spectral sequence converses to 
$HH^*(C^*(M); C^*(N))$ as an algebra. 
Moreover 
the K{\"u}nneth theorem yields that the $E_1$-term is a complex of the form
$$
E_1^{p,q}=\text{Hom}_{A\otimes A^{op}}^{p+q}
(\B_p(H^*(M); H^*(M); H^*(M)), H^*(N))
$$ which is the Hochschild cochain complex.
We have the result.
\end{proof}

The spectral sequence in Theorem \ref{thm:Mss} 
is called {\it the  Moore spectral sequence}.

Let ${\bf C}^*$
be the Hochschild complex mentioned in the proof of 
Theorem \ref{thm:Mss}.   
Then the inclusion $i : F^p{\bf C}^n \to {\bf C}^n$ defines the submodule
$$
F^pHH^n:= \text{Im}\{H(i) : H^n(F^p{\bf C}^*) \to HH^n(C^*(M);C^*(N))\}
$$
of the Hochschild cohomology $HH^n(C^*(M);C^*(N))$. We define  
{\it the associated bigraded module} of $HH^n(C^*(M);C^*(N))$ by
$$
\text{Gr}^{p,q}HH^*(C^*(M); C^*(N)) = F^pHH^{p+q}/F^{p+1}HH^{p+q}.
$$
Observe that $\text{Gr}^{p,q}HH^*(C^*(M); C^*(N))$ is isomorphic to
the vector space $E_\infty^{p, q}$ in the $E_\infty$-term 
of the Moore spectral sequence.

We prove a key proposition to proving the assertion in Theorem
\ref{thm:2} concerning the bigraded algebra structure. 

\begin{prop}
\label{prop:key}Let $M$ be a simply-connected Poincar{\'e} duality space
and $\{\widetilde{E}_r^{*,*}, \widetilde{d}_r\}$ the Eilenberg-Moore spectral sequence associated
with the pull-back diagram
$$
\xymatrix@C25pt@R25pt{
LM \ar[r] \ar[d] & M^{[0, 1]} \ar[d]^{\e_0\times \e_1} \\
M \ar[r]^{\Delta} & M\times M
}
$$
converging to the cohomology $H^*(LM; \K)$, where $\e_i$ denotes the
evaluation map at $i$ for $i = 0, 1$ and $\Delta$ is the diagonal map.
Let $\{E_r, d_r\}$ be the Moore spectral sequence converging to
$HH^*(C^*(M); C^*(M))$. 
Then all the elements in the $E_2$-term of 
$\{\widetilde{E}_r^{*,*}, \widetilde{d}_r\}$ with total degree less than
or equal to $l$ are permanent cycles if and only if so are all the
elements in the $E_2$-term of $\{E_r, d_r\}$ with total degree greater
than or equal to $-l+d-1$. 
In particular, the spectral sequence 
$\{\widetilde{E}_r^{*,*}, \widetilde{d}_r\}$ 
collapses at the $E_2$-term if and only if
so does the Moore spectral sequence $\{E_r, d_r\}$.
\end{prop}

\begin{proof}
Let $\{\B_p, d\}_{p \geq 0}$ be the normalized bar complex of $A:=H^*(M)$.
The $E_1$-term of the spectral sequence 
$\{\widetilde{E}_r^{*,*}, \widetilde{d}_r\}$ is given by
$\widetilde{E}_1^{-p, q} = (\B_p\otimes_{A\otimes A^{op}}A)^{-p+q}$.
Thus we see that the spectral sequence 
$\{(\widetilde{E}_r^{*,*})^\vee, \widetilde{d}_r^\vee\}$, 
which is dual to $\{\widetilde{E}_r^{*,*}, \widetilde{d}_r\}$, 
converges to $HH^*(C^*(M); C^*(M)^\vee)$ and that 
$$
(\widetilde{E}_1^{-p, q})^\vee = \text{Hom}^{p-q}(\B_p\otimes_{A\otimes A^{op}}A,
\K)
\cong \text{Hom}^{p-q}_{A\otimes A^{op}}(\B_p, A^\vee) 
$$
as complexes. The main theorem  in \cite{J} asserts that as a vector space
$H_{-p+q}(LM)\cong HH^{p-q}(C^*(M), C^*(M)^\vee)$. 

Moreover, since $A$ is a commutative, it follows that 
the Poincar\'e duality gives an isomorphism 
$A \stackrel{\cong}{\to} A^\vee$ of $A$-bimodules. Then we have an isomorphism 
\begin{eqnarray*}
(\widetilde{E}_2^{-p, q})^\vee = H(\widetilde{E}_1^{-p, q})^\vee
&=&H(\text{Hom}^{p-q}_{A\otimes A^{op}}(\B_p, A^\vee)) \\
&\cong&
H(\text{Hom}^{p-q+d}_{A\otimes A^{op}}(\B_p, A)) = H(E_1^{p,-q+d})
=E_2^{p,-q+d}.
\end{eqnarray*}
The result \cite[13 Theorem]{F-T-V} due to F\'elix, Thomas and Vigu\'e-Poirrier allows us
to obtain an isomorphism
$HH^*(C^*(M); C^*(M)^\vee) \cong HH^{*+d}(C^*(M); C^*(M))$;
see also \cite[Theorem 20]{Luc2}. Thus it
turns out that the following conditions are equivalent:

\smallskip
\noindent
(1) all the elements in the $E_2$-term of 
$\{\widetilde{E}_r^{*,*}, \widetilde{d}_r\}$ with total degree less than
or equal to $l$ are permanent cycles. \\
(2) $H^n(LM)\cong \displaystyle\bigoplus_{-p+q=n} 
\widetilde{E}_2^{-p, q}$ for any
$n\leq l$. \\
(3) $H_n(LM)\cong
\displaystyle\bigoplus_p H^{-n}(\text{Hom}_{A\otimes A^{op}}(\B_p,
A^\vee))$
for any $n\leq l$. \\
(4) $HH^{-n+d}(C^*(M); C^*(M))\cong HH^{-n}(C^*(M); C^*(M)^\vee)
\cong \displaystyle\bigoplus_{p-q=-n} E_2^{p, -q+d}$ for any $n\leq l$. \\
(5) all the
elements in the $E_2$-term of $\{E_r, d_r\}$ with total degree greater
than or equal to $-l+d-1$ are permanent cycles. 

\smallskip 
\noindent
In fact, the isomorphisms mentioned above allow us to conclude that 
the conditions (2), (3) and (4) are equivalent.   
Since the vector space $E_{r+1}^{p,q}$ is a subquotient of $E_r^{p,q}$,
 more precisely, 
$$E_{r+1}^{p,q}\cong 
\text{Ker} \{d_r : E_r^{p,q}\to E_r^{p+r, q+1-r}\}/
\text{Im} \{d_r : E_r^{p-r, q-1+r}\to E_r^{p, q}\},  
$$ 
it follows that  
$\dim E_{r+1}^{p,q} \leq \dim E_r^{p,q}$. Then the equality holds if and
 only if 
all the elements in $E_r^{p,q}$ and $E_r^{p-r, q-1+r}$ are cycles. 
This yields that $\dim E_\infty^{p,q} = \dim E_2^{p,q}$ for 
$p+q \leq l$ if and only if the all the elements in $E_2^{p,q}$ are
 permanent cycles for  $p+q \leq l$. The fact implies that 
the conditions (1) and (2) are equivalent. The same argument does work
 well to show the equivalence of (4) and (5).       
We have the result. 
\end{proof}

\begin{prop}
\label{prop:3}
Let $X$ be a simply-connected space as in Theorem \ref{thm:2}. Then
$$
\text{\em Gr}HH^*(C^*(X;{\mathbb Z}/2); C^*(X;{\mathbb Z}/2)) \cong
\wedge(y_1, y_2, ..., y_l)\otimes {\mathbb Z}/2[\nu_1^*, \nu_2^*, ...,
\nu_l^*]
$$
as a bigraded algebra, where $\text{\em bideg} \  y_j =(0, \deg y_j)$ and 
$\text{\em bideg} \ \nu_j^* =(1,  - \deg y_j )$.
\end{prop}

\begin{proof}
Since $Sq^1\equiv 0$ on $H^*(X; {\mathbb Z}/2)$ by assumption, 
it follows from 
\cite[Theorem]{Smith2} that the Eilenberg-Moore
spectral sequence converging to
$H^*(LX; {\mathbb Z}/2)$ collapses at the $E_2$-term. Thanks to
Proposition \ref{prop:key}, we see that
the Moore spectral sequence for $C^*(X)$ collapses at the
$E_2$-term. Proposition \ref{prop:HH} yields
that $E_\infty^{*, *} $ is isomorphic to $H^*(X)\otimes {\mathbb
Z}/2[\nu_1^*, ...., \nu_l^*]$ as a bigraded
algebra. This completes the proof.
\end{proof}

In the case where the characteristic of the underlying field is odd, 
we can solve the extension problem in the associated bigraded algebra
for an appropriate space.

\begin{prop}
\label{prop:4}
Let $p$ be odd prime and 
$G$ a simply-connected H-space whose mod $p$ cohomology is an
exterior algebra, 
say $H^*(G; {\mathbb Z}/p)\cong \Lambda (y_1, ..., y_l)$. Then
$$
HH^*(C^*(G;{\mathbb Z}/p); C^*(G;{\mathbb Z}/p)) \cong
\wedge(y_1, y_2, ..., y_l)\otimes {\mathbb Z}/p[\nu_1^*, \nu_2^*, ...,
\nu_l^*]
$$
as an algebra, where $\text{\em bideg} \  y_j =(0, \deg y_j)$ and 
$\text{\em bideg} \ \nu_j^* =(1,  - \deg y_j )$.
\end{prop}

\begin{proof}
The spectral sequence $\{\widetilde{E}_r^{*,*}, \widetilde{d}_r\}$ converging to 
$H^*(LG; {\mathbb Z}/p)$ collapses at the $E_2$-term. In fact, since $G$
is an H-space, it follows that $LG$ is homotopy equivalent to the
product $G\times \Omega G$, where $\Omega G$ denotes the based loop
space. This implies that 
\begin{eqnarray*}
H^*(LG; {\mathbb Z}/p)&\cong& H^*(G;{\mathbb Z}/p)\otimes 
H^*(\Omega G;{\mathbb Z}/p) \\
&\cong& \wedge(y_1, y_2, ..., y_l)\otimes 
\Gamma [\nu_1, ..., \nu_l] \cong \text{Total}\widetilde{E}_2^{*,*}, 
\end{eqnarray*}
where $\deg \nu_i = \deg y_i-1$. 
The third isomorphism follows from the usual computation of 
the $E_2$-term with the Koszul-Tate resolution described  
in Proposition \ref{prop:resolution}. 
By virtue of Proposition \ref{prop:key}, we see that the Moore spectral sequence 
also collapses at the $E_2$-term. Thus Proposition \ref{prop:HH} yields that 
$$
\text{Gr}HH^*(C^*(G;{\mathbb Z}/p); C^*(G;{\mathbb Z}/p)) \cong
\wedge(y_1, y_2, ..., y_l)\otimes {\mathbb Z}/p[\nu_1^*, \nu_2^*, ...,
\nu_l^*]
$$
as a bigraded algebra. It is immediate that $y_i^2=0$ 
in $HH^*(C^*(G;{\mathbb Z}/p); C^*(G;{\mathbb Z}/p))$  
for any $i$ because $p$ is odd. 
We have the result. 
\end{proof}

\section{The associated bigraded Batalin-Vilkovisky algebra}
In this section, a Batalin-Vilkovisky algebra structure on
the Moore spectral sequence are considered. We first recall  
the Batalin-Vilkovisky algebra structure on the Hochschild cohomology 
defined in \cite{Luc2}. 

Let $M$ be a simply-connected Poincar\'e duality space of formal
dimension $d$ and let $A$ stand for the singular cochain algebra
$C^*(M;\K)$. Let ${\mathbb B}$ denote the normalized bar complex 
${\mathbb B}(A;A;A)$.  
We define an isomorphism of complexes 
$$
\iota : \text{Hom}(A\otimes_{A\otimes A^{op}}{\mathbb B}, \K) \stackrel{\cong}{\to}
\text{Hom}_{A\otimes A^{op}}({\mathbb B}, A^\vee)
$$
by $\iota(f)(\alpha)(a) = (-1)^{|a||\alpha |}f(a\otimes \alpha)$ for 
$\alpha \in  {\mathbb B}$ and $a \in A$. Here the $A$-bimodule structure 
of $A^\vee$ is defined by 
$
\langle f\cdot \alpha \cdot g ; h\rangle=
(-1)^{|f|}\langle \alpha; ghf\rangle
$ 
for $f, g, h \in A$ and $\alpha \in A^\vee$. 
Then one obtains an isomorphism 
$$
\xymatrix@C20pt@R15pt{
\iota^* : \text{Hom}(H(A\otimes_{A\otimes A^{op}}{\mathbb B}), \K)) &
H(\text{Hom}(A\otimes_{A\otimes A^{op}}{\mathbb B}, \K)) \ar[r]^(0.6){H(\iota)} _(0.6){\cong}
\ar[l]_(0.45){\kappa}^(0.45){\cong} & HH^*(A; A^\vee),  
}
$$
where $\kappa$ denotes the K\"unneth isomorphism. 
Observe that the source of the map $\iota^*$ is the dual 
$HH_*(A; A)^\vee$ to the Hochschild homology $HH_*(A; A)$ of $A$.  
We also recall the quasi-isomorphism 
$J : A\otimes_{A\otimes A^{op}}{\mathbb B} \to C^*(LM)$
of differential graded modules due to Jones \cite{J}. Then it follows  
that this quasi-isomorphism fits in the commutative diagram
$$
\xymatrix@C20pt@R15pt{C^*(LM) & & A\otimes_{A\otimes A^{op}}{\mathbb B} \ar[ll]_{J}^{\simeq} \\
& C^*(M) \ar[lu]^{ev^*} \ar[ru]_{\eta'}, 
}
$$ 
where $ev : LX \to X$ is the evaluation map at zero and $\eta'$ is the chain map defined by 
$\eta'(a)=a\otimes 1$. Therefore we have a commutative diagram 
$$
\xymatrix@C20pt@R25pt{
H^*(LM)^\vee  \ar[r]^{H(J)^\vee}_{\cong} \ar[rd]_{H(ev)^\vee} & HH_*(A;A)^\vee \ar[d]_{H(\eta')^\vee}& 
H(\text{Hom}(A\otimes_{A\otimes A^{op}}{\mathbb B}, \K)) \ar[r]^(0.6){H(\iota)}_(0.6){\cong} \ar[l]_(0.6){\kappa} ^(0.6){\cong}
\ar[d]^{H({\eta'}^\vee)} & HH^*(A;A^\vee) \ar[d]^{HH(\eta, 1)}\\
& H^*(M)^\vee & H(\text{Hom}(A, \K)) \ar[l]_{\kappa}^{\cong} \ar[r]^(0.5){H(\iota)}_(0.5){\cong} &HH^*(\K; A^\vee), 
}
\eqnlabel{add-0}
$$
where $\eta : \K \to A$ denotes the unit. It is readily seen that 
a section $s : M \to LM$ 
of the evaluation map $ev$ induces a section $H(s)^\vee$ of the map 
$H(ev)^\vee$. 
Let 
$B$ be the Connes boundary map on 
$A\otimes T(s\bar{A})\cong A\otimes_{A\otimes A^{op}}{\mathbb B}$; see \cite{G-J}. 
By definition, we see that  
$$
B(a_0[a_1| a_2| ... |a_k]) =\sum_{i=0}^{k}(-1)^{(\e_i+1)(\e_{k+1}-\e_i)}1[a_i| ... | a_k | a_0 | ... |a_{i-1}].
$$
We then have

\begin{prop} \label{prop:Luc2}
\cite[Propositions 11 and 12]{Luc2}{\em (i)} Let $\omega_A^\vee \in H(A)^\vee$ be the dual base of the 
fundamental class of $M$. Define an element
$[m] \in HH^{-d}(A, A^\vee)$  
by $[m]=\iota^* H(J)^\vee H(s)^\vee(\omega_A^{\vee})$. Then the product $- \smile [m]$
induces an isomorphism
$$
\theta : HH^p(A; A) \to HH^{p-d}(A; A^\vee).
$$
{\em (ii)} The Hochschild cohomology ring  $HH^*(A; A)$ is a Batalin-Vilkovisky
algebra equipped with the B-V operator
$\Delta$ of degree $-1$ defined by the composite
$$
\xymatrix@C25pt@R25pt{
HH^{p}(A;A) \ar[r]^(0.45){\theta}_(0.45){\cong} &
HH^{p-d}(A; A^\vee) \ar[r]^{{\iota^*}^{-1}} &
HH_{-p+d}(A;A)^\vee \ar[d]_{{H(B)}^\vee} \\
HH^{p-1}(A;A) \ar[r]_(0.45){\theta}^(0.45){\cong}
& HH^{p-d-1}(A; A^\vee) &
HH_{-p+d+1}(A;A)^\vee \ar[l]^{\iota^*} . 
}
$$
\end{prop}

\begin{rem}
\label{rem:PD}
Let $H$ be a Poincar\'e duality algebra with the fundamental class $\omega_H$ 
and $PD :  H \stackrel{\cong}{\to} H^\vee$ the isomorphism 
of $H$-bimodules defined by the Poincar\'e duality; that is, 
$PD(1)=\omega_H^\vee$, where $\omega_H^\vee$ denotes 
the dual element to  $\omega_H$. 
As mentioned in the proof of Proposition \ref{prop:key}, 
the map $PD$ induces the isomorphism 
$HH^*(1; PD) : HH^*(H; H) \stackrel{\cong}{\to} HH^*(H; H^\vee)$. Moreover, 
it is readily seen  that the isomorphism $HH(1; PD)$ coincides 
with the cup product 
$$
\smile \omega_H^\vee : HH^*(H; H) \to HH^*(H, H^\vee)
$$   
Here $\omega_H^\vee$ is considered an element in 
$H^\vee \cong HH^{0,*}(H; H^\vee)$.  
\end{rem}

We retain the same notations as in the proof 
of Theorem \ref{thm:Mss}. 
Let $\widehat{\bf C}$ and $\overline{\bf C}$ stand for the Hochschild cochains 
$\text{Hom}_{A\otimes A^{op}}({\mathbb B}, A^\vee)$ and
$\text{Hom}_{\K\otimes \K}(\K, A^\vee)$, respectively. 
We define the same filtrations $\{F^p\widehat{\bf C}\}$ 
and $\{F^p\overline{\bf C}\}$ 
as that of ${\bf C}$. 
Their filtrations  construct the spectral sequences 
$\{\widehat{E}_r^{*,*}, \widehat{d}_r\}$ and 
$\{\overline{E}_r^{*,*}, \overline{d}_r\}$ 
converging to $HH^*(A; A^\vee)$ and 
$HH^*(\K, A^\vee)=H(A^\vee)$, respectively. 
We see that $F^p\overline{\bf C}=0$ for $p>0$ and hence 
$\overline{E}_r^{p,*}=0$ for $p>0$. 
Since the cup product respects the filtrations 
$\widehat{\bf C}$ and $\overline{\bf C}$, 
it follows that the spectral sequences converge to the targets as algebras.  
Indeed, the target $A^\vee$ of the Hochschild cochains should be
replace with $(C^*(M)/I)^\vee$ as in the proof of 
Theorem \ref{thm:Mss} when considering the convergence of the spectral
sequence. However to simplify, we also write $A^\vee$ 
for the reduction. 

We give a B-V algebra structure to the Moore spectral sequence.          

\begin{thm} 
\label{thm:MSS} Let $M$ be a simply-connected Poincar\'e duality space. Then
the Moore spectral sequence $\{E_r^{*,*}, d_r\}$ converging to
$HH(C^*(M); C^*(M))$ admits the structure of 
a differential Batalin-Vilkovisky bigraded algebra, 
in the sense that each term $E_r^{*,*}$ is endowed with 
the B-V operator $\Delta _r: E_r^{p, q} \to E_r^{p-1, q}$ such that
$d_r\Delta_r + \Delta_r d_r = 0$, $H(\Delta_r)=\Delta_{r+1}$ and
$E_\infty^{*,*}$ is isomorphic to $\text{\em Gr}HH^*(C^*(M); C^*(M))$
as bigraded Batalin-Vilkovisky algebras. 
\end{thm}

\begin{proof} 
We first recall that the $E_r$-term of the Moore spectral sequence
 $\{E_r, d_r\}$ is defined by 
$E_r^{s, t}=Z_r^{s, t}/(Z_{r-1}^{s+1, t-1}+B_{r-1}^{s, t})$, where 
$Z_r=F^s{\bf C}^{s+t}\cap \partial^{-1}(F^{s+r}{\bf C}^{s+r+1})$ 
and 
$B_{r}^{s, t}=F^s{\bf C}^{s+t}\cap \partial(F^{s-r}{\bf C}^{s+r-1})$. 
The $E_r$-term $\widehat{E}_r^{s, t}$ of the spectral sequence 
$\{\widehat{E}_r^{*,*}, \widehat{d}_r\}$
is defined by the same form with 
the filtration $\{F^p\widehat{\bf C}\}$.  

Let $m \in F^0\widehat{\bf C}^{-d}$ 
be a cocycle representing the element $[m] \in HH^{-d}(A, A^\vee)$ 
described in Proposition \ref{prop:Luc2}. 
Then it follows from \cite[Lemma 2.1]{K-S} that  
$\{m\}$ is a permanent cycle. The cup product 
$$
\ \smile \ : \text{Hom}_{A\otimes A^{op}}(\B_*, A)
\otimes \text{Hom}_{A\otimes A^{op}}(\B_*, A^\vee) \to 
\text{Hom}_{A\otimes A^{op}}(\B_*, A^\vee)
$$
respects the filtrations; that is, 
$F^s{\bf C}^n \smile F^t\widehat{\bf C}^m \subset 
F^{s+t}\widehat{\bf C}^{n+m}$. 
Therefore the product with the element 
$\{m\}\in \widehat{E}_2^{0, -d}\cong \widehat{E}_r^{0, -d}$ 
induces a morphism 
$$
E(m)_r:= - \smile \{m\} : E_r^{p,q} \to \widehat{E}_r^{p,q-d}
$$  
of spectral sequences. We show that $E(m)_2$ is an isomorphism 
and hence so is $E(m)_r$ for $2 \leq r \leq \infty$. 

We observe that the unit $\eta : \K \to A$ induces the morphism 
$\{E(\eta)_r\} : \{\widehat{E}_r^{*,*}, \widehat{d}_r\} 
\to \{\overline{E}_r^{*,*}, \overline{d}_r\}$ of spectral sequences.  
Consider the map 
$$
E(\eta)_2 : \widehat{E}_2^{0,-d}\cong 
HH^{0,-d}(H(A), H(A^\vee)) \to \overline{E}_2^{0,-d}\cong
H^{-d}(A^\vee). 
$$ 
Then it follows that the K\"unneth map 
$\kappa : H(\text{Hom}(A; \K)) \stackrel{\cong}{\to} \text{Hom}(H(A); \K) = H(A)^\vee$ 
sends the image $E(\eta)_2(\{m\})$ of $\{m\} \in \widehat{E}_2^{0,-d}$ 
to the dual to the fundamental class.  
To see this, we consider the commutative diagram
$$
\xymatrix@C30pt@R17pt{
HH^{-d}(A; A^\vee) \ar@{->>}[d]_\pi \ar[r]^(0.55){HH(\eta; 1)} & 
H^{-d}(A^\vee) \ar@{=}[d] \\
\widehat{E}_\infty^{0,-d} \ar@{->}[d]_i \ar[r]^{E(\eta)_\infty} & 
\overline{E}_\infty^{0,-d} \ar@{=}[d] \\
\widehat{E}_2^{0,-d} \ar[r]_{E(\eta)_2} & \overline{E}_2^{0,-d},
}
\eqnlabel{add-1}
$$
where $\pi$ and $i$ denote the natural projection and the natural inclusion, 
respectively.  
Then the definition of $[m]$ and the diagram (4.1) enable us to deduce
that 
$$
\kappa E(\eta)_2(\{m\})= \kappa E(\eta)_2i\pi([m]) 
= \kappa HH(\eta; 1)[m] = \omega_M^\vee.
$$
Moreover, we see that 
the map $E(\eta)_2: H^*(A^\vee)\cong \widehat{E}_2^{0, *} \to 
\overline{E}_2^{0,*}=H^*(A^\vee)$ 
is the identity since $E(\eta)_2=HH(H(\eta), 1)$.
This implies that $\kappa (\{m\}) =\kappa E(\eta)_2(\{m\})= \omega_M^\vee$. 
Therefore, the map $E(m)_2$ coincides with the cup product with 
$\omega_M^\vee$. Remark \ref{rem:PD} 
allows one to conclude that $E(m)_2$ is an isomorphism. 

Define a map $\Delta' : \widehat{\bf C} \to \widehat{\bf C}$ 
with degree $-1$ by $\Delta' 
= \iota\circ B^\vee \circ \iota^{-1}$, where 
$B^\vee$ is defined by $B^\vee(f)=(-1)^{|f|}f\circ B$. 
It is readily seen that $\Delta'$ maps 
$F^p\widehat{\bf C}^n$ to $F^{p-1}\widehat{\bf C}^{n-1}$ and hence 
the map $\Delta'$ induces a morphism 
$\{E(\Delta')_r\} : \{\widehat{E}_r^{*,*}, \widehat{d}_r\} 
\to \{\widehat{E}_r^{*,*}, \widehat{d}_r\}$
of spectral sequences with bidegree $(-1, 0)$.  Thus 
we define maps $\Delta_r : E_r^{*, *} \to E_r^{*, *} $ by 
$\Delta_r = E(m)_r^{-1}\circ \Delta_r' \circ E(m)_r$ for $r \geq 2$, 
which give the morphism of spectral sequences from 
$\{E_r^{*,*}, d_r\}$ to itself with bidegree $(-1, 0)$. 
Proposition \ref{prop:Luc2}(ii) implies that 
the $E_2$-term of the Moore spectral sequence admits the structure of
a Batalin-Vilkovisky algebra  
with the operator $\Delta_2$. 
The equality $bB+Bb =0$ for the differential $b$ of 
the Hochschild complex $A\otimes T(s\bar{A})$ enables us to deduce 
that  $d_r\Delta_r + \Delta_r d_r = 0$ and that
$H(\Delta_r)=\Delta_{r+1}$. 

The isomorphism 
$E_\infty^{*,*} \stackrel{\cong}{\to} \text{Gr}HH^*(C^*(M); C^*(M))$ 
is induced by the inclusion of cocycles. 
Therefore the isomorphism is compatible with the B-V operators. 
Observe that the filtration $\{F^pHH^n\}$ of $HH^n(A; A)$ and the
filtration $\{F^p\widehat{H}^n\}$ of $HH^n(A; A^\vee)$ are  bounded for
each $n$ since $H^*(M)$ is of finite dimension; see the proof of 
Theorem \ref{thm:Mss}.  
Recall that
$E(m)_\infty : E_\infty^{*,*} \to \widehat{E}_\infty^{*,*}$ 
is an isomorphism. 
This yields that the map from $\{F^pHH^n\}$ to $\{F^p\widehat{H}^n\}$ 
induced by $\theta$  
is also an isomorphism. This completes the proof.
\end{proof}

We have fundamental properties of the Moore spectal sequence. 

\begin{cor} 
\label{cor:4.4}
With the same notations as in the proof of Theorem \ref{thm:MSS}, 
one obtains that \\
{\em (i)} the map 
$$
E(\eta)_2 : \widehat{E}_2^{0,*}\cong HH^{0,*}(H(A), H(A^\vee))=H^*(A^\vee)
\to \overline{E}_2^{0,*}\cong H^*(A^\vee)
$$ 
is the identity and the element $\{m\}\in \widehat{E}_2^{0, -d}$
is the dual to the fundamental class $\omega_M$ in $H(A)=H^*(M; \K)$, 
and that \\
{\em  (ii)} each element in $E_2^{0,*}$ is a permanent cycle. 
\end{cor}

\begin{proof} 
The assertion (i) follows from the argument with the diagram (4.2) 
in the proof of Theorem \ref{thm:MSS}. 

We prove the assertion (ii). 
It suffices to show that all the elements in $\widehat{E}_2^{0,*}$ 
are permanent cycles because  $E(m)_r$ is an isomorphism for any $r\geq 2$.
The result follows from the commutative diagram
$$
\xymatrix@C25pt@R20pt{
\widehat{E}_2^{0,*} \ar[d]_{E(\eta')_2}^\cong 
& \  \widehat{E}_3^{0,*} \ar@{>->}[l]  &
\ \cdots \ar@{>->}[l]  & \ \widehat{E}_\infty^{0, *} \ar@{>->}[l]  & 
HH^*(A, A^\vee) . \ar@{->>}[l] \ar@/^1pc/[dllll]^{HH^*(\eta, 1)} \\
H^*(A^\vee)
}
$$
In fact, since the evaluation map $LM \to M$ admits a section, 
it follows from the diagram  (4.1) that 
the map $HH^*(\eta, 1)$ is an epimorphism. We have the result. 
\end{proof}

In order to prove Theorem \ref{thm:2}, 
we look at the behavior of the  B-V operator on the $E_2$-term of 
the Moore spectral sequence. 
To this end, we give an explicit isomorphism between 
the Hochschild homology of a graded module computed by the bar resolution and 
that computed by the Koszul-Tate resolution described in Proposition \ref{prop:resolution}. 

\begin{lem}
\label{lem:bar-KT}
Let $A$ be an exterior algebra over a field of characteristic $2$, say $A\cong \wedge(y_1, ..., y_l)$. 
Then there exists an isomorphism 
$$
\phi : H(A\otimes T(s\bar{A}), b) \to 
H(A\otimes_{A\otimes A}{\mathcal F}, 1\otimes d)
\cong \wedge(y_1, ..., y_l)\otimes \Gamma[\nu_1, ..., \nu_l]
$$
such that $\phi([y_i])=\nu_i$, $\phi([y_i|y_i]) = \gamma_2(\nu_i)$ and 
$\phi([y_i|y_j]+[y_j|y_i])=\nu_i\nu_j$ for $i\neq j$. 
Here $(A\otimes T(s\bar{A}), b)$ and $({\mathcal F}, d)$ 
denote the Hochschild complex and the resolution of $A$ 
mentioned in Proposition \ref{prop:resolution}, respectively.  
\end{lem}

\begin{proof}
We construct a DG $A\otimes A^{op}$-module map $\xi=\{\xi_i\}:
{\mathbb B}(A; A ;A) 
\to {\mathcal F}$ which covers the identity map on $A$ by induction 
on the filtration degree. 

Let ${\mathcal F}_i $ denote a submodule of the form 
$\{ x \in {\mathcal F} \ | \  \text{bideg} \ x = (-i, *) \}$. 
Suppose that a map $u : Z_i \to {\mathcal F}_i$ satisfies the condition that 
$\partial u = \xi_{i-1}\partial$, where $Z_i$ is an 
$A\otimes A^{op}$-module which is a direct summand of $\B_i$; that is, 
$\B_i= Z_i\oplus Z_i'$ for some $A\otimes A^{op}$-submodule $Z_i'$. 
Since $\B_i$ is a free $A\otimes A^{op}$-module and ${\mathcal F}$ 
is acyclic, it follows that there exists an $A\otimes A^{op}$-module map 
$\xi_i : \B_i \to {\mathcal F}_i $ such that 
$\xi_i|_{Z_i} = u$ and $\partial \xi_i = \xi_{i-1}\partial$. 
Observe that $Z_i'$ is projective.

Define $\xi_0 : A\otimes A ={\mathbb B}_0(A; A; A) 
\to {\mathcal F}_0=A\otimes A^{op}$ 
to be the identity map and 
$\xi_1 :  {\mathbb B}_1(A; A; A)=A\otimes T^1(s\bar{A})\otimes A 
\to A\otimes A\otimes \K\{\nu_i ; i = 1, .., l\}$ by extending 
the $A\otimes A$-module map which sends basis $[y_i]$ to $\nu_i$.
We observe that 
$d\xi_1([y_i])=d\nu_i = y_i\otimes 1 + 1\otimes y_i=
\xi_0d([y_i])$. 
A direct computation shows that  
\begin{eqnarray*}
\xi_1d([y_i|y_j]+[y_j|y_i])&=&(y_i\otimes 1 + 1\otimes y_i)\nu_j +
(y_j\otimes 1 + 1\otimes y_j)\nu_i  \ \text{for} \ i\neq j \ \ \text{and} \\
\xi_1d([y_i|y_i])&=&(y_i\otimes 1 + 1\otimes y_i)\nu_i.
\end{eqnarray*}
This enables us to define an $A\otimes A^{op}$-module map 
$\xi_2 : {\mathbb B}_2(A; A; A)\to {\mathcal F}_2$ which is compatible
with the differentials and satisfy the condition that   
$\xi_2([y_i|y_i])=\gamma_2(\nu_i)$ and
$\xi_2([y_i|y_j]+[y_j|y_i])=\nu_i\nu_j$. We define 
$\xi : {\mathbb B}(A; A ;A) 
\to {\mathcal F}$ extending $\{\xi_i\}_{0\leq i \leq 2}$ to 
the whole complex ${\mathbb B}(A; A ;A)$. 

Let $\phi' : A\otimes T(s\bar{A}) \to A_{A\otimes A^{op}}{\mathbb B}(A; A ;A)$ be 
an isomorphism of complexes defined by 
$\phi'(a_0[a_1| ...|a_n])=a_0\otimes 1[a_1| ...|a_n]1$. The induced map
$\phi=H((1\otimes \xi) \circ \phi') : H(A\otimes T(s\bar{A}), b) 
\to H(A\otimes_{A\otimes A}{\mathcal F}, 1\otimes d)$ 
is the desired isomorphism. 
\end{proof}

\noindent
{\it Proof of Theorem \ref{thm:2}.} 
It follows from the proof of Proposition \ref{prop:3} that 
as bigraded algebras  
$$
\text{Gr}HH^*(C^*(X;{\mathbb Z}/2); C^*(X;{\mathbb Z}/2)) \cong
E_2^{*,*} \cong 
\wedge(y_1, y_2, ..., y_l)\otimes {\mathbb Z}/2[\nu_1^*, \nu_2^*, ...,
\nu_l^*],
$$ where $\text{bideg} \  y_j =(0, \deg y_j)$ and 
$\text{bideg} \ \nu_j^* =(1,  - \deg y_j )$.
Thus in order to prove Theorem \ref{thm:2}, it suffices to determine 
the B-V structure on the $E_2$-term of the Moore spectral sequence.  

Let $H$ denote the cohomology $H^*(M)$ and let $\omega_H$ be 
the fundamental class, namely $\omega_H=y_1\cdots y_l$. 
We first recall the isomorphism of complexes 
$$
\iota : \text{Hom}(H\otimes_{H\otimes H^{op}}\F, \K) \to 
\text{Hom}_{H\otimes H^{op}}(\F, H^\vee)
$$
defined by $\iota(f)(\alpha)(a)=(-1)^{|a||\alpha|}f(a\otimes \alpha)$
for $a \in H$ and $\alpha \in \F$, where $\F \to H \to 0$ is the
Koszul-Tate resolution of $H$ mentioned in Proposition
\ref{prop:resolution}. Recall also the $H$-bimodule structure on
$H^\vee$; see the beginning of this section. 
Using the isomorphism $\iota$, the cup product $- \smile \omega_H^\vee$
and $\phi$ in Lemma \ref{lem:bar-KT}, we can determine the B-V operator
on the $E_2$-term. 
The key to the computation is that the Connes boundary map $B$ is a
derivation on cycles modulo boundary with respect to the shuffle product on 
the Hochschild complex. 
By virtue of Corollary \ref{cor:4.4} (i), 
we have 
$\theta(\nu_i^*\nu_j^*)=\nu_i^*\nu_j^*\smile {m}=\nu_i^*\nu_j^*\smile \omega_H^\vee
= \omega_H^\vee\nu_i^*\nu_j^*$, which sends 
$\nu_i\nu_j$ to $\omega_H^\vee \in H^\vee$. 
Moreover, by the definition of $\iota$, one obtains 
$$
\iota((\omega_H\nu_i\nu_j)^*)(\nu_i\nu_j)(a)=
(\omega_H\nu_i\nu_j)^*(a\nu_i\nu_j) = \left\{
\begin{array}{ll}
1 & \ \text{if} \ a = \omega_H \\
0 & \ \text{otherwise}. 
\end{array}
\right.
$$
This yields that
$\theta(\nu_i^*\nu_j^*)=\iota(\omega_H\nu_i\nu_j)^*$. 
We then have 
\begin{eqnarray*}
& & \langle H(B)^\vee \phi^\vee (\omega_H\nu_i\nu_j)^*, 
y_{i_1}\cdots y_{i_u}[y_t] \rangle \\
&=&\langle (\omega_H\nu_i\nu_j)^*, 
\phi H(B)(y_{i_1}\cdots y_{i_u}[y_t]) \rangle \\
&=&\langle (\omega_H\nu_i\nu_j)^*, 
\phi \sum_{i_s}y_{i_1}\cdots \check{y}_{i_s}\cdots
y_{i_u}[y_{i_s}]*[y_t]
\rangle \\
&=& 
\langle (\omega_H\nu_i\nu_j)^*, 
\sum_{i_s\neq t}y_{i_1}\cdots \check{y}_{i_s}\cdots
y_{i_u}\phi([y_{i_s}|y_t]+[y_t|y_{i_s}])
\rangle \\
&=& 
\langle (\omega_H\nu_i\nu_j)^*, 
\sum_{i_s\neq t}y_{i_1}\cdots \check{y}_{i_s}\cdots
y_{i_u}\nu_{i_s}\nu_t
\rangle = 0.  
\end{eqnarray*}
Here $*$ denotes the shuffle product on the Hochschild homology and 
$\check{y}_j$ means that the element $y_j$ has been deleted.   
The second equality follows from the fact that $H(B)$ is a derivation 
with respect to the shuffle product; see \cite[Lemma 4.3]{G-J}. 
In fact, we see that 
\begin{eqnarray*}
H(B)(y_{i_1}\cdots y_{i_u}[y_t]))
&=&H(B)(y_{i_1}* \cdots * y_{i_u}*[y_t]) \\
&=&(H(B)y_{i_1})*y_{i_2}*\cdots * y_{i_u}*[y_t] \\
&& + \cdots + y_{i_1}*\cdots *y_{i_{u-1}}* (H(B)y_{i_u})*[y_t] \\ 
&& + \ y_{i_1}*\cdots * y_{i_u}* (H(B)[y_t]) \\
&=& 
\sum_{i_s} y_{i_1}* \cdots *[y_{i_s}]* \cdots * y_{i_u}*[y_t] + 
y_{i_1}* \cdots * y_{i_u}*0 \\
&=&
\sum_{i_s}y_{i_1}\cdots \check{y}_{i_s}\cdots
y_{i_u}[y_{i_s}]*[y_t]. 
\end{eqnarray*}
Thus we have $\Delta(\nu_i^*\nu_j^*)=0$. 

We compute  $\Delta(y_j\nu_i^*)$, $\Delta(\nu_i^*)$ and 
$\Delta(y_i)$ below.  Since 
$\omega_H^\vee \cdot y_j = (y_1\cdots \check{y}_j \cdots y_l)^*$  
in $H^\vee$, it follows that 
\begin{eqnarray*}
\theta(y_j\nu_i^*)&=& \omega_H^\vee (y_j\nu_i^*) \\
&=& (y_1\cdots \check{y}_j \cdots y_l)^*\nu_i^* \\
&=& \iota ((y_1\cdots \check{y}_j \cdots y_l \nu_i)^*).  
\end{eqnarray*}
Moreover, we see that 
\begin{eqnarray*}
& & \langle H(B)^\vee \phi^\vee ((y_1\cdots \check{y}_j \cdots y_l\nu_i)^*), 
\omega_H \rangle \\
&=& \langle (y_1\cdots \check{y}_j \cdots y_l \nu_i)^*, 
\phi \sum_{t=1}^l(y_1\cdots \check{y}_t \cdots y_l)[y_t]\rangle 
= \delta_{ij}.  
\end{eqnarray*}
The fact that $\theta(1)=\omega_H^\vee$ allows us to deduce that 
$\Delta(y_j\nu_i^*)=\delta_{ij}\cdot 1$. 
For dimensional reasons, we have $\Delta(\nu_i^*)=\Delta(y_j)=\Delta(y_iy_j)=0$. 
This completes the proof. 
\hfill\qed


\medskip

The following corollary illustrates that the Moore spectral sequence is
reliable when calculating explicitly 
the Hochschild cohomology of the singular cochain on a space.  

\begin{cor}
\label{cor:S-times-S}
Let $M$ be a simply-connected mod $2$ Poincar\'e duality space 
whose mod $2$ cohomology is isomorphic to 
an exterior algebra of the form $\wedge(y_1, y_2)$, 
where $\deg y_1= \deg y_2 = n$. Suppose that $n > 4$. Then 
as a Batalin-Vilkovisky algebra
$$
HH^*(C^*(M;{\mathbb Z}/2); C^*(M;{\mathbb Z}/2)) \cong
\wedge(y_1, y_2)\otimes {\mathbb Z}/2[\nu_1^*, \nu_2^*]
$$
in which $\Delta(y_j)=0$, $\Delta(y_jy_j)=0$, $\Delta(\nu_i^*)=0$,
$\Delta(\nu_i^*\nu_j^*)=0$ 
for $1\leq i, j \leq l$ and 
$\Delta(y_i\nu_j^*)=\delta_{ij}\cdot 1$, where 
$\deg y_j =n$ and $\deg   \nu_j^* =-n +1$ for $j = 1$ and $2$.
\end{cor}


\begin{proof}
By virtue of Theorem \ref{thm:2}, we see that
as bigraded Batalin-Vilkovisky algebras  
$$
\text{Gr}HH^*(C^*(M;{\mathbb Z}/2); C^*(M;{\mathbb Z}/2)) \cong
E_\infty^{*,*} \cong 
\wedge(y_1, y_2)\otimes {\mathbb Z}/2[\nu_1^*, \nu_2^*] 
$$
with $\Delta(y_j)=0$, $\Delta(\nu_i^*)=0$, 
$\Delta(\nu_i^*\nu_j^*)=0$ for $1\leq i, j \leq l$ and 
$\Delta(y_i\nu_j^*)=\delta_{ij}\cdot 1$. Observe that 
$Sq^1=0$ on $H^*(M; {\mathbb Z}/2)=\wedge (y_1, y_2)$. 

We have to solve extension problems on the product and on the B-V operator. 
Since there exists no nonzero element in $E_\infty^{p,q}$ for $p\geq 1$
and $p + q = 2n$, it follows that $y_i^2=0$ for $i=1$ and $2$; see
the figure displayed below. 
$$
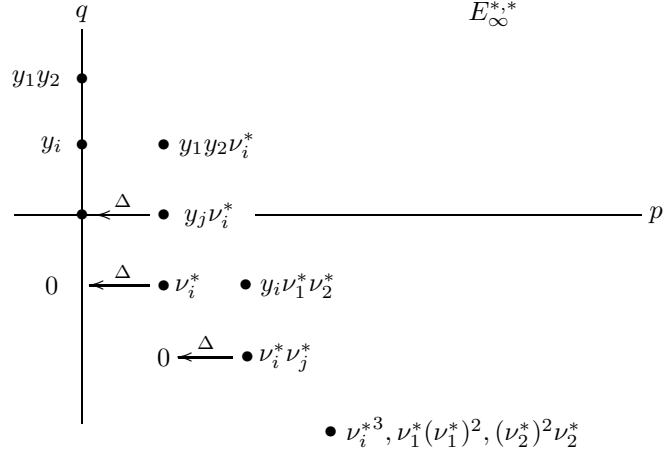
\xymatrix@C20pt@R11pt{
{} & {} q \ar@{-}[dddddd] & &  & & E^{*,*}_\infty \\
 y_1y_2 \hspace*{-0.8cm}& {\bullet}  \ar@{-}[ddd] & {}  & \\
 y_i  \hspace*{-1.2cm}  & \bullet &  \bullet  
                                 & \hspace*{-0.8cm}{y_1y_2\nu_i^*} &  \\
\ar@{-}[rr] &  \bullet  & \bullet\ar[l]_{\Delta}  
                                 & \hspace*{-1.2cm}{y_j\nu_i^*} 
                     \ar@{-}[rrr] \hspace*{-0.2cm} &  &  &  p& \\
0 \hspace*{-1.2cm} &     & \bullet \ar[l]_{\Delta} 
                                 & \hspace*{-0.9cm}{\nu_i^*}
                                     \hspace*{0.5cm} \bullet  &   
                      & \hspace*{-5.2cm}{y_i\nu_1^*\nu_2^*} &   \\
 &      {} & 0          & \bullet \ar[l]_{\Delta} 
                                 & \hspace*{-1.2cm} \nu_i^*\nu_j^*   \\
& {} &   &    &  \bullet & \hspace*{-0.8cm}{{\nu_i^*}^3}, 
\nu_1^*(\nu_1^*)^2, (\nu_2^*)^2\nu_2^* & 
}
$$

We consider the extension problems on the B-V operator. 
For non-positive integers $l$ and $m$ with $l+m\geq 3$, we see that 
\begin{eqnarray*}
\deg y_i{\nu_1^*}^l{\nu_2^*}^m < \deg y_1y_2{\nu_1^*}^l{\nu_2^*}^m 
&=& 2n+(l+m)(-n+1)  \\
&=& 2+(l+m-2)(-n+1) \leq 2+ (-n+1) < 0,  
\end{eqnarray*}
$\deg y_i{\nu_1^*}{\nu_2^*}=-n +2$, $\deg y_1y_2{\nu_1^*}{\nu_2^*}=2$, 
$\deg y_i{\nu_j^*}=1$ and  $\deg y_1y_2{\nu_i^*}=n+1$.    
This enables us to conclude that 
$\Delta(y_i\nu_j^*)=\delta_{ij}\cdot 1$, 
$\Delta(y_j)=0$ and 
$\Delta(y_iy_j)=0$
on $HH^*(C^*(M); C^*(M))$ 
because $\deg \Delta(y_iy_j)=2n-1$, $\deg \Delta(y_i)= n-1$ and $\deg\Delta(y_i\nu_j^*)=0$. 
For $l\geq 3$, one has 
\begin{eqnarray*}
\deg y_1y_2{\nu_i^*}^l-(\deg {\nu_1^*}{\nu_2^*} - 1)
&=& 2n+l(-n+1) -2(-n+1)+1 \\
&=&\left\{
\begin{array}{ll}
n+2 & \ \text{if} \ l=3 \\
3   & \ \text{if} \ l=4 \\ 
\text{a negative integer} & \  \text{if} \ l>4
\end{array}
\right. .
\end{eqnarray*}
Moreover, we see that 
\begin{eqnarray*}
\deg y_i{\nu_j^*}^l-(\deg {\nu_1^*} {\nu_2^*} - 1)
&=& 2n+l(-n+1) -2 (-n+1)+1 -n \\
&=&\left\{
\begin{array}{ll}
2 & \ \text{if} \ l=3 \\
3-n   & \ \text{if} \ l=4 \\ 
\text{a negative integer} & \  \text{if} \ l>4. 
\end{array}
\right. 
\end{eqnarray*}
Since $E_\infty^{2, q} =0$ for $q< -2n$ and 
$\Delta({\nu_1^*}{\nu_2^*}) =0$ 
in $E_\infty^{1, -2n} $, 
it follows that $\Delta({\nu_1^*}{\nu_2^*}) =0$ on $HH^*(C^*(M); C^*(M))$. 
The same calculation as above shows that    
\begin{eqnarray*}
\deg y_1y_2{\nu_i^*}^l-(\deg {\nu_i^*}  - 1)
&=&\left\{
\begin{array}{ll}
n+2 & \ \text{if} \ l=2 \\
3   & \ \text{if} \ l=3 \\ 
\text{a negative integer} & \  \text{if} \ l>3. 
\end{array}
\right. 
\end{eqnarray*}
and that
\begin{eqnarray*}
\deg y_i{\nu_j^*}^l-(\deg {\nu_i^*}  - 1)
&=&\left\{
\begin{array}{ll}
2 & \ \text{if} \ l=2 \\
3-n   & \ \text{if} \ l=3 \\ 
\text{a negative integer} & \  \text{if} \ l>3. 
\end{array}
\right. 
\end{eqnarray*}
These deduce that $\Delta({\nu_i^*})=0$  on $HH^*(C^*(M); C^*(M))$. 
We have the result. 
\end{proof}

\begin{rem}
In the case $n=3$, we write $\Delta(y_j\nu_i^*) 
= \delta_{ij}\cdot 1 + \e y_1y_2{\nu_1^*}^l{\nu_2^*}^m$, where 
$l+m=3$. The filtration argument we use above does not 
work well to determine whether $\e$ is zero or not. 
\end{rem}

\noindent
{\it Acknowledgments}.
The author thanks Hiroshi Nagase for useful comments on 
the ordinary Hochschild cohomology rings of a polynomial algebra and 
an exterior algebra. 
He also thanks the referee for careful reading of 
a previous version of this paper and 
for valuable suggestions to revise the article.

\section{Appendix}

In this short section, 
we summarize the notion of ``a ring 
homomorphism'' from the Hochschild cohomology to the graded center of a
triangulated category.  Though there is {\it no result} yet, we have   
a certain expectation that string topology plays a crucial role in the study 
of global structure of the cochain algebra on a space via the homomorphism.

We begin by recalling the definition of the graded center of 
a triangulated category. 

\begin{defn} (cf.\cite[3.2]{B-F}, \cite[\S 2]{Linckel})
Let ${\mathcal T}$ be a $\K$-linear triangulated category with
suspension functor $\Sigma$.  
The graded center ${\mathcal Z}({\mathcal T})$ is a graded  
family whose degree $n$ component ${\mathcal Z}^n({\mathcal T})$
consists of all natural transformations 
$\varphi : Id_{\mathcal T} \to \Sigma^n$ such that 
$\varphi \Sigma = (-1)^n \Sigma \varphi$. 
\end{defn} 

Let $R$ be a commutative graded ring and 
$\Phi : R \to {\mathcal Z}({\mathcal T})$ a ring homomorphism preserving 
the degree.  
Here we ignore set theoretic issues on the graded center. Indeed, 
the ring homomorphism 
means that, for each object $X$ in ${\mathcal T}$, one has a
homomorphism of graded algebra 
$\Phi_X : R \to \text{End}_{\mathcal T}^*(X)$ 
such that 
$$\Phi_Y(\alpha)\beta = (-1)^{|\alpha||\beta|}\beta \Phi_X(\alpha)
$$ 
for $\alpha \in R$ and $\beta \in \text{Hom}_{\mathcal T}^*(X, Y)$.  

Let $A$ be a DG algebra over a field $\K$. 
Then we have a triangulated category $\D(A)$, which is
the derived category of DG modules over $A$ with the shift functor $\Sigma$; 
$(\Sigma N)^n = N^{n+1}$, as the suspension functor. 
It follows from  \cite[Proposition 1.1]{BGSS} 
that the cup product on $HH^*(A; A)$ 
coincides with the Yoneda product. 
We then have a ring homomorphism $\Phi$ from the Hochschild cohomology
ring $HH^*(A; A)$ 
to the graded center 
of the triangulated category $\D(A)$. 
In fact, the homomorphism $\Phi : HH^*(A; A) \to {\mathcal Z}(\D(A))$ 
is defined by 
$$
\Phi(f)(M)=\Phi_M(f)=Id_M\otimes_A f : M \to \Sigma^nM
$$ 
in $\D(A)$ for $f \in HH^n(A, A)$.

Let $X$ be a simply-connected space whose cohomology with coefficients in
$\K$ is locally finite. The general argument 
above gives a ring homomorphism 
$$
\Phi : 
HH^*(C_*(\Omega X), C_*(\Omega X)) \to  
{\mathcal Z}(\D(C_*(\Omega X))). 
$$
By virtue of 
\cite[Theorem I]{F-H-T_COBAR}, we have a quasi-isomorphism 
of DG algebras 
from the cobar complex of $C_*(X)$ to $C_*(\Omega X)$.   
Thus the result \cite[Theorem 1]{F-M-T} due to F\'elix, Menichi and Thomas
allows one to obtain an isomorphism of algebras 
between $HH^*(C_*(\Omega X), C_*(\Omega X))$
and $HH^*(C^*(X), C^*(X))$, which indeed respects 
the Gerstanhaber algebra structure. Suppose further that $X$ is a closed
oriented manifold. 
Then we have a ring homomorphism form the loop
homology ${\mathbb H}_*(LX)$ to the graded center of the derived
category $\D(C_*(\Omega X))$ with the isomorphism between the loop
homology and the Hochschild homology \cite{C, C-J, F-T-V2, Me}.  

The same way allows us to define a ring homomorphism   
$$
{\mathbb H}_*(LX)\cong HH^*(C^*(X), C^*(X)) \to {\mathcal Z}(\D(C^*(X))). 
$$
Thus it is expected that the loop homology 
is of great use when studying triangulated categories
associated with cochain algebras on spaces via the theory of support
varieties; see for example \cite{B-I-K_L, B-I-K, S-S}, and when
considering the (co)chain type level \cite{K2, K3} of a space 
over a simply-connected manifold. Topological properties of spaces,
which the graded center captures, remains to be elucidated.

\end{document}